\documentclass[12pt,a4paper,leqno]{amsart}
\usepackage{amsmath,amsfonts}
\usepackage{graphicx}

\newtheorem*{theoremintrod}{Theorem 7.1}

\newtheorem{theorem}{Theorem}[section]
\newtheorem{proposition}[theorem]{Proposition}

\newtheorem{lemma}[theorem]{Lemma}
\newtheorem{corollary}[theorem]{Corollary}
\newtheorem{definition}{Definition}[section]

\usepackage{amssymb}
\usepackage{upref}
\usepackage[mathcal]{eucal}

\title{Homogeneous Weyl connections of non-positive curvature}
\author{Gabriela Tereszkiewicz and Maciej P. Wojtkowski}
\address{Faculty of Mathematics and Computer Science
\newline University of Warmia and Mazury
\newline 10-710 Olsztyn, POLAND}
\thanks{The authors were  partially supported 
by the NCN grant 2011/03/B/ST1/04427} 
\email{wojtkowski@matman.uwm.edu.pl}
\date{\today} 
\subjclass{53C24 (53C21)}

\begin{document}

\begin{abstract}
We study homogenous Weyl connections with non-positive 
sectional curvatures. The Cartesian product $\mathbb S^1 \times M$
carries  canonical families of Weyl connections with such a property, 
for any Riemmanian manifold $M$. 
We prove that if a homogenous Weyl connection on a manifold, modeled
on a unimodular Lie group, is non-positive in a stronger sense
(streched non-positive), then it must be locally of the product type. 
\end{abstract}

\maketitle

\section{Introduction}
Let $M$ be a Riemannian manifold with the metric $g$. 
A Weyl connection on $M$ is
a linear symmetric connection with conformal parallel transport. 
A modern account of the theory of such connections was given by 
Folland, \cite{F}.

Let us recall that a Weyl 
connection $\widehat\nabla$  is determined by the 1-form $\varphi$
satisfying $\widehat{\nabla}_X g = -2\varphi(X)g$. The difference 
of two linear connections is a tensor.  Denoting by $\nabla$
the Levi-Civita connection, we get 
\[
\widehat{\nabla}_XY = \nabla_XY + \varphi(Y)X + \varphi(X)Y - 
\langle X , Y\rangle E, 
\]
for any tangent vector fields $X,Y$, where $E$ is the vector field
dual to $\varphi$, and $\langle\cdot,\cdot\rangle$ denotes the metric
$g$.

For any smooth function $U$ on $M$ the pair $(g,\varphi) $ defines 
the same Weyl connection as the pair  $(e^{-2U}g, \varphi+dU)$.
In particular for a closed 1-form $\varphi$ the Weyl connection  
is locally a Levi-Civita connection of some metric in the conformal class.

The curvature operator of a Weyl connection  
$$\widehat{R}(X,Y) = \widehat{\nabla}_X\widehat{\nabla}_Y - 
\widehat{\nabla}_Y\widehat{\nabla}_X -
\widehat{\nabla}_{[X,Y]}$$ has the symmetric part,
called the distance curvature,
$\widehat{R}_s(X,Y) = d\varphi(X,Y) I$, where $I$ is 
the identity operator, \cite{F}. 
The anti-symmetric part of the curvature operator 
$\widehat{R}_a = \widehat{R} - \widehat{R}_s 
$, called the direction curvature, can be used
to define the sectional curvature  $\widehat K(\Pi)$
of the Weyl connection in the direction of a plane 
$\Pi$
$$
\widehat K(\Pi) = \langle \widehat{R}_a(X,Y)Y,X\rangle,
$$
for any orthonormal basis $\{X,Y\}$ of $\Pi$.
These sectional curvatures depend 
on the choice of the Riemannian metric $g$ in the conformal 
class, and not on the Weyl connection alone.
Indeed, although the anti-symmetric part of the curvature operator 
$\widehat{R}_a$ is independent of the choice of the metric,
with the passage to the metric $e^{-2U}g$ the sectional curvatures
acquire the universal factor $e^{2U}$.
Hence the sign of sectional curvatures is well defined.

For a Weyl connection on a compact manifold Gauduchon, \cite{Gau}, 
established the existence and
uniqueness of a metric in the conformal class such that 
the  respective vector field $E$ is divergence free, $div E =0$. 
We will say in such a case that the pair $(g,E)$ is
the {\it Gauduchon representation} of the Weyl connection, 
and $g$ is the {\it Gauduchon metric}, 
irrespective of the manifold $M$ being compact or not.

It was observed in \cite{W1,W2} that Weyl connections are related
to a special class of dynamical systems called Gaussian thermostats.
The paper of Gallavotti and Ruelle, \cite{G-R}, gives
an account of the physical significance of these systems.
A Gaussian thermostat is defined  by a Riemannian metric on a manifold $M$
and a vector field $E$. The trajectories of the Gaussian thermostat
satisfy the following ordinary differential equations in the tangent
bundle $TM$,\cite{W1}. 
\begin{equation}
\frac{d}{dt}x = v, 
\nabla_vv = E-\frac{\langle E,v\rangle}{\langle v,v\rangle}v,
\end{equation}
where $x = x(t)\in M$ is a parametrized curve in $M$.
By the force of these equations the ``kinetic energy'' $v^2$ is constant. 
Fixing the value of this constant $k = v^2$ we obtain that the trajectories
of (1) are geodesics of the Weyl connection defined by the vector field
$k^{-1}E$. (Note that this coincidence does not extend to the parametrizations:
in general the natural parameter of the Weyl geodesic is not
proportional to the arc length.)

We can see that the behavior of the trajectories of a Gaussian 
thermostat for different values of the kinetic energy is governed by the 
family of Weyl connections defined by the family of vector fields
$\gamma E, \gamma > 0,$ where large  values of $\gamma$ correspond to 
small values of the kinetic energy.

The goemetry and dynamics of multidiemsional Gaussian thermostats were studied
recently in the paper of Daraibekov and Paternain, \cite{D-P}.

In the present paper we study Weyl connections with non-positive sectional 
curvatures, the {\it non-positive} Weyl connections. 
The motivation comes from the results of 
\cite{W1,W2,P-W}. 
In the Riemannian case if all the sectional curvatures are 
negative then the  geodesic flow is an Anosov system. 
It was proven
in \cite{W1} that negative sectional curvatures of the Weyl connection lead
to some hyperbolic properties, namely the dominated splitting with exponential 
growth and decay of volumes. A simplified account of this result 
can be found in \cite{P-W}.      

In particular in dimension 2 we get the Anosov property.
In this case there is only one  sectional curvature
and it is equal to the Gaussian curvature of the Gauduchon metric. 
Surfaces with negative Gaussian curvature are well understood. For any
divergence free vector field $E$ on such a surface
we get the Gaussian thermostat with the Anosov property.

In dimension 3 and higher we do not know of examples of Weyl connections
with negative sectional curvatures,
beyond the small perturbations of Levi-Civita connections with 
negative sectional curvatures.

In the search  of such examples one encounters rigidity phenomena
similar to the Riemannian case.

In the work of Alekseevski\cite{Al}, Heintze \cite{He},  and
Azencott and Wilson \cite{Az-W}, the homogeneous Riemannian metrics 
with non-positive sectional curvatures were classified. In particular 
it was found that left invariant metrics on unimodular Lie groups 
are by necessity flat.

Constant vector fields on flat tori give examples of non-positive 
Weyl connections. It follows from the results in \cite{W3} that 
there are no other non-positive Weyl connections on tori with 
the flat Gauduchon metric. Moreover there are partial results in
\cite{W3} towards the  conjecture that also for  
arbitrary Gauduchon metrics
there are no other non-positive Weyl connections on tori. 
For example it was proven in \cite{W2} 
that there are no locally Riemannian non-positive Weyl connections
on tori except for the ones defined by constant vector fields.  
Weyl connections on $\mathbb S^1\times \mathbb S^{n-1}$ were also studied 
in \cite{W3}. In the present paper we generalize these results.

The examples in \cite{W3} have the additional properties that 
the vector field $E$ is {\it parallel}, i.e., $\nabla E = 0$,
and it is {\it stretched non-positive}
\[
\begin{aligned}
(SNP) \ \ &\text{there is} \ \ \gamma_0 \geq 0 \ \ \text{such that for every} 
\ \ \gamma \geq \gamma_0  \\ 
&\text{the Weyl connection defined by} 
\ \ \gamma E \ \ \text{is non-positive}.
\end{aligned}
\]

It is straightforward that any parallel vector field $E$ is  
stretched non-positive (cf. Proposition 2.3).

Note that the stretched non-positivity of a vector field $E$
is a property of the Gaussian thermostat rather than the respective 
Weyl connection. Indeed for two different representations  
$(g_1,E_1), (g_2,E_2)$ of the same Weyl connection the families  
$\gamma E_1$ and $\gamma E_2$ 
define in general different families of Weyl connections. 
We do not know under what conditions the two families share 
non-positivity for large values of $\gamma$. 
We will not address such issues in this paper. 

Nevertheless, in the case of a compact manifold we define 
{\it a stretched non-positive Weyl connection} 
by the requirement that in the Gauduchon 
representation $(g,E)$ ($div \ E =0$), the field $E$ is SNP.

Let us note here the paper of Narita, \cite{N}, where 
the Gauduchon representation is used to define Weyl connections
with constant sectional curvatures.

In Section 2 we consider Weyl connections defined by Killing 
vector fields, especially the non-positive connections. 
Parallel vector fields play an important role in this section.
As a bi-product of this discussion we obtain a notable restriction 
on the Clifford isometries on manifolds with some non-positive 
sectional curvature (Proposition 2.2).

In Section 3 we study the notion of stretched non-positivity
(SNP) for vector fields, and for Weyl connections, on general 
Riemannian manifolds.

In Section 4 homogeneous Weyl connections on Lie groups are 
introduced, with low dimensional examples. We exhibit 
 an isolated non-positive Weyl connection on SOL,\cite{T},
and SNP Weyl connections on groups with the euclidean metric (unimodular),
and with the hyperbolic metric (nonunimodular).

Section 5 contains 
results about multidimensional examples of solvable extensions 
of abelian groups. We obtain there families of examples of 
isolated non-positive Weyl connections.
It would be interesting to classify this type of 
geometric objects.

In Section 7 we present our main result which says that
the only homogeneous SNP vector fields (or Weyl connections)
modeled on unimodular Lie groups are parallel.

Our result should be compared with the 
Riemannian case where, on one hand, there are
no non-positive, non-flat, 
Levi-Civita connections on unimodular Lie groups.
On the other hand there are many homogeneous 
metrics with negative sectional curvature
on manifolds modeled on unimodular Lie groups.

In Section 6 we introduce the homogeneous Weyl connections of a fairly 
general type. They were considered previously in the paper by Kerr, \cite{Ke}. 
We will follow here the notations of the book by Besse, \cite{Be}. 
Let $G$ be a Lie group and $H$ its compact subgroup. Let 
$\mathfrak g = \mathfrak h \oplus \mathfrak 
p$ be an $Ad(H)$ invariant
splitting of the Lie algebra $\mathfrak g$ of $G$, into the Lie algebra
$\mathfrak h$ of $H$ and a complementary linear subspace $\mathfrak p$. 
Note that $\mathfrak p$ is not  a Lie subalgebra in general. 
Any choice of an $Ad(H)$ invariant scalar product in $\mathfrak p$ 
gives rise to a $G$ invariant metric on the homogeneous
space $M = G/H$. We consider the subspace
\[
\mathfrak p_0 = \{X \in \mathfrak p| [Y,X] = 0 \ \ \forall Y\in \mathfrak h\},
\]
on which the action of $Ad(H)$ is trivial.
We choose a left invariant vector field $E$ defined by an element 
from $\mathfrak p_0$.
It projects onto $M = G/H$ as a $G$-invariant vector field. The subspace
$\mathfrak p_0$ may be equal to zero, however 
for $\mathfrak h = 0$ we have that $\mathfrak p_0 =\mathfrak g$. 
The last case is that of a left invariant metric on $G$, 
with a left invariant vector field $E$.

Our main result says

\begin{theoremintrod}
If the Lie group $G$ is unimodular and a $G$-invariant vector 
field on $M= G/H$, defined by $E \in \mathfrak p_0$, is SNP
then the vector field is parallel.
\end{theoremintrod}

We include in Section 6 a derivation of the 
formula for the sectional curvature of 
a homogeneous Riemannian metric, using the Jacobi
fields. This derivation clarifies the calculations
appearing in the proof of the main result in Section 7.

Throughout the paper we will assume that the dimension of $M$ is at least $3$.
In dimension 2 for a Gauduchon representation the Weyl connection 
is non-positive if and only if the Gauss curvature of the metric 
is non-positive (see the formula (4) for the sectional curvature 
of the Weyl connection in Section 2).

\section{Killing fields and non-positive Weyl connections }\label{ssnp}

Let us consider a Killing field $E$ on a Riemannian manifold $M$,
i.e. for every vector fields $X$ and $Y$ 
\[
\left< \nabla_XE,Y\right> + \left< \nabla_YE,X\right> = 0.
\]

\begin{proposition}\label{wektKill}
For a  Killing vector field $E$ the following are equivalent

(i) $E$ has constant length,

(ii) $\nabla_EE = 0$,

(iii) integral curves of $E$ are geodesics.
\end{proposition}
\begin{Proof}
If $E$ is a Killing vector field then for any  vector field $X$ we have
\begin{equation}
\nabla_XE^2 = 2\left< \nabla_XE,E\right> = 
-2\left< \nabla_EE,X\right>.
\end{equation}
It shows the equivalence of (i) and (ii).

If (iii) holds then we have for some positive smooth function $f$
$\nabla_E(fE) = 0$. It follows that $\nabla_E E$ is parallel to $E$.
However for a Killing vector field $\nabla_E E$ must be orthogonal to $E$. 
Hence (ii) follows.
\end{Proof}$\square$

Killing vector fields with constant length define one parameter
groups of Clifford isometries, i.e. isometries which move all points
by the same distance. 
It is well known that there are no Clifford 
isometries on manifolds of negative curvature (reference?). We have further

\begin{proposition}\label{prodKill}
Let $E$ be a Killing vector field of constant length on a 
Riemannian manifold $M$. If at every point of $M$ the sectional curvatures
of all planes containing $E$ are non-positive, then all these 
sectional curvatures must be actually equal to zero and
for every vector field $X$,  $\nabla_X E = 0$.
\end{proposition}
\begin{Proof}
Let us assume for simplicity that $ E^2 =1$. Denote by 
$\psi^t, t\in \mathbb R,$ the one parameter
group of isometries generated by $E$. 
By Proposition~\ref{wektKill} the integral curves of $E$ are 
geodesics and further $\psi^t$ shifts points by the distance $t$ along these
geodesics. Let us consider a family of such geodesics $\gamma(t,u)$
where $t\in \mathbb R$ is the arc length parameter along a geodesic, 
and $u \in (-\epsilon,\epsilon)$ is the parameter of the family. 
We get the Jacobi field 
$
J(t) =\frac{\partial}{\partial u} \gamma(t,u)_{|u=0}.
$
Since $\psi^{t_1}\gamma(t_2,u) = \gamma(t_1+t_2,u)$
we get the invariance of the Jacobi field under the flow
$\psi^t$, i.e., 
$
D\psi^tJ(0) = J(t).
$
Hence the Jacobi field has constant length.

On the other hand it follows from the Jacobi equations that 
\begin{equation}
\nabla_EJ = Y = \nabla_JE \ \ \ 
\text{and} 
\ \ \ 
\frac{d^2}{dt^2}\frac{1}{2}J^2 = Y^2 - k J^2,
\end{equation}
where $k$ denotes the sectional curvature of the plane spanned by
$E$ and $J$. Since the Jacobi field has constant length the last 
expression must be equal to zero, i.e.,
$
 Y^2 = k J^2.
$ 
Under the assumption $k\leq 0$ we conclude that $Y = 0$ and $k=0$
everywhere. 
\end{Proof}$\square$

We say that a vector field $E$ is {\it parallel} if 
for every vector field $X$,  $\nabla_XE = 0$. It is well known
that the presence of a parallel field forces the local Cartesian
product structure of $M$, see \cite{K-N}, Chapter 99.
Because of the crucial role it plays 
in our discussion we include a simple direct 
argument. 
\begin{proposition}\label{parKill}
If a vector field $E$ is parallel then
its orthogonal distribution is integrable and
the Riemannian manifold has locally the Cartesian product 
structure of an orthogonal section of $E$ by an integral curve
of $E$. In particular the curvature operator $R(X,E) =0$ for any
vector field $X$. 
\end{proposition}
\begin{Proof}
By the assumption the field $E$ is parallel along any curve.
Hence the parallel transport along any curve preserves 
also the orthogonal complement of $E$.
It follows that for $X$ and $Y$,  two vector fields orthogonal to $E$,
the derivative   $\nabla_XY$ is orthogonal to $E$. Finally 
the Lie bracket $[X,Y] =  \nabla_XY -\nabla_YX$ is orthogonal 
to $E$. By the Frobenius Theorem the orthogonal distribution is
integrable. 

Let us consider a local orthogonal section $N$ of the flow of $E$,
i.e., a local integral submanifold of the orthogonal distribution.
Since the isometries in the flow translate points from $N$ along
the geodesics by a constant distance, it follows that locally
the Remannian metric is isometric to that of the Cartesian product
$(-\epsilon,\epsilon)\times N$. 

The last claim follows immediately from the product structure.
\end{Proof}$\square$

It is clear that conversely for the Cartesian product 
$N\times \mathbb R$ the ``vertical'' vector field 
of constant length is  parallel.

As described in the Introduction, we will  call a Weyl
connection {\it non-positive} if it has only non-positive
sectional curvatures.

For a Weyl connection defined by a vector field $E$
the sectional curvature $\widehat K= \widehat K(\Pi)$
in the direction of a plane $\Pi$ is equal to (\cite{W1,P-W})
\begin{equation}
\widehat{K}(\Pi) = K(\Pi) - E^2_{\perp} - \ div_{\Pi}E,
\end{equation}
where $E_{\perp}$ is the component of $E$ orthogonal to the plane
$\Pi$ and the {\it partial divergence} $div_{\Pi}E$ is equal to
\[
div_{\Pi}E = \langle \nabla_XE,X\rangle +\langle \nabla_YE,Y\rangle,
\]
for any orthonormal frame $(X,Y)$ of the plane $\Pi$.
Moreover $K(\Pi)$ denotes the Riemannian sectional curvature in the plane
$\Pi$.

Let us assume that the field $E$ is a Killing vector field.
In such a case the formula (4) reads
$
\widehat{K}(\Pi) = K(\Pi) - E^2_{\perp},
$
and for planes $\Pi \ni E$ we have
$
\widehat{K}(\Pi) =  K(\Pi).
$
Hence to get a non-positive Weyl connection defined by a Killing
vector field $E$ we have to have $K(\Pi) \leq 0$ for every plane 
$\Pi \ni E$. If we assume further that $E$ has constant
length we conclude by Proposition~\ref{parKill} that 
$E$ is parallel and the Riemannian manifold has locally 
the Cartesian product structure. Conversely, parallel vector fields 
always have the {\it SNP property}. Namely we have the following
\begin{proposition}\label{snp}
For a unit parallel vector field $E$ the Weyl connection 
defined by $\gamma E$ is non-positive for 

$\gamma^2 \geq
\max\{ K(\Pi)|\ \Pi \ \ \text{orthogonal to} \ \  E\}$.
\end{proposition}
We will be using several times the following formulas
\begin{lemma}
Let 
$X, Y_1,Y_2$ be mutually orthogonal (locally defined)
unit vector fields on $M$. 
For a plane $\Pi$ spanned by $aX+bY_1, a^2+b^2=1$
and $Y_2$, the Riemannian sectional curvature $K(\Pi)$
in the direction of $\Pi$ is equal to
\[
\begin{aligned}
&K(\Pi) = \left<R(Y_2, aX+bY_1)(aX+bY_1),Y_2\right> \\ 
=&a^2\left<R(Y_2, X)X,Y_2\right> +
b^2\left<R(Y_2, Y_1)Y_1,Y_2\right>
+ 2ab\left<R(Y_2, X)Y_1,Y_2\right>.
\end{aligned}
\]
Further for  a smooth function $f$ and the vector field $E = fX$ 
its partial divergence
$div_{\Pi} E$ in the direction of $\Pi$,
is equal to
\[
\begin{aligned}
&div_{\Pi}E =
 \left<\nabla_{aX+bY_1}(E), aX+bY_1 \right>  +
\left<\nabla_{Y_2}(E), Y_2 \right> =
\\
&a^2(\left<\nabla_{X}E,X \right>  + \left<\nabla_{Y_2}E, Y_2 \right>)+
b^2(\left<\nabla_{Y_1}E, Y_1 \right> + \left<\nabla_{Y_2}E, Y_2 \right>)\\
&+ab(\left<\nabla_{X}E,Y_1 \right>+\left<\nabla_{Y_1}E, X \right>)=\\
&a^2\left(df(X) + f\left<\nabla_{Y_2}X, Y_2 \right>\right)+
b^2\left(f\left<\nabla_{Y_1}X, Y_1 \right> + 
f\left<\nabla_{Y_2}X, Y_2 \right>\right)\\
&+ab\left(f\left<\nabla_{X}X,Y_1 \right>+df(Y_1)\right)
\end{aligned}
\]
\end{lemma}$\square$

\begin{Proof} (Proposition 2.4)
Let us consider three mutually orthogonal unit vectors
$E, Y_1, Y_2$. We consider the
Riemannian sectional curvature $K(\Pi)$ in the direction of the plane 
$\Pi$ spanned by $Y_2$ and $aE+bY_1, a^2 + b^2 =1$.
Since for any vector field $Y$ the curvature operator $R(Y,E) = 0$
we have by Lemma 2.5
\begin{equation}
\label{ks}
\begin{aligned}
K(\Pi) =
&a^2\left<R(Y_2, E)E,Y_2\right> +
b^2\left<R(Y_2, Y_1)Y_1,Y_2\right>
+\\ &ab(\left<R(Y_2, E)Y_1,Y_2\right> + 
\left<R(Y_2, Y_1)E,Y_2\right>\\
 =&b^2\left<R(Y_2, Y_1)Y_1,Y_2\right>.
\end{aligned}
\end{equation}
Hence the curvature of the Weyl connection defined by $\gamma E$ is 
equal to 
\[
\widehat{K}(\Pi) = b^2\left<R(Y_2, Y_1)Y_1,Y_2\right>
-b^2\gamma^2.
\]
\end{Proof}$\square$

Let us now assume that a compact manifold $M$ possesses 
a parallel vector field $E$. We will look for non-positive 
Weyl connections other than the ones described in 
Proposition~\ref{snp}. First of all if we multiply 
$E$ by a function $f$ then we loose the non-positivity 
unless the function is constant.

\begin{proposition}
\label{pnp}
For a unit parallel vector field $E$ on a compact manifold
and a non-constant smooth function $f$ 
the Weyl connection defined by $ fE$ has some positive sectional 
curvatures.
\end{proposition}
\begin{Proof}
As in the proof of Proposition~\ref{snp} we calculate the sectional 
curvature of the plane $\Pi$ spanned by $aE+bY_1$ and  $Y_2$.
We have further by Lemma 2.5
$
div_{\Pi}(fE) = a^2\nabla_Ef +ab\nabla_{Y_1}f.
$
Hence 
\[
\widehat{K}(\Pi) = b^2(\left<R(Y_2, Y_1)Y_1,Y_2\right>
-f^2) - a^2\nabla_Ef - ab\nabla_{Y_1}f.
\]
If the connection is non-positive then $\nabla_Ef$ is non-negative
everywhere. However for any vector field  $E$ with zero divergence
 on a compact manifold and any smooth function $f$ the integral of 
$\nabla_Ef$ is zero. It follows that if the connection is non-positive
then $\nabla_Ef$ is equal to zero.

Now non-positivity forces also $\nabla_{Y_1}f =0$ everywhere,
for any vector field $Y_1 \perp E$. Hence the function $f$ is
constant.
\end{Proof}$\square$

For a Riemannian manifold with a  unit parallel field $E$ 
the Ricci curvature $Ric(aE+bY) = b^2 Ric(Y)$, 
for any unit vector field $Y\perp E$, $a^2+b^2 =1$.
If we assume that the Ricci curvature is positive except for $E$
then the only non-positive Weyl connections are those described in
Proposition~\ref{snp}.

\begin{proposition}\label{fsnp}
If for a compact manifold $M$ with a unit parallel vector field $E$ 
the Ricci curvature $Ric(Y) > 0$ for every unit vector $Y\perp E$
then the only non-positive Weyl connections on $M$ are those defined by 
$\gamma E$, where $\gamma^2 \geq
\max_{\Pi \perp E} K(\Pi)$.
\end{proposition}
\begin{Proof}
Let $F$ be a vector field defining a non-positive Weyl connection.
For a plane $\Pi$ spanned by $F$ and a unit vector $Y\perp F$
\[
\widehat{K}(\Pi) = K(\Pi) -
\left<\nabla_XF, X\right> - \left<\nabla_YF, Y\right> \leq 0,
\]
where $X = \frac{1}{|F|}F$. Adding the above inequalities
over an orthogonal basis $(Y_1,Y_2,\dots, Y_{n-1},F)$ we obtain 
\[
Ricc(X) - div F - (n-2)\left<\nabla_XF, X\right> \leq 0.
\]
This inequality is meaningful at every point where $F \neq 0$.
Multiplying the inequality by $|F|^{n-2}$ we arrive at 
\[
|F|^{n-2}Ricc(X) - |F|^{n-2} div F - 
(n-2)|F|^{n-2}\left<\nabla_XF, X\right> \leq 0,
\]
which holds everywhere on $M$, with the first term vanishing
where $F =0$.

Let us introduce the vector field $V =  |F|^{n-2}F$. For $n\geq 3$
the field $V$ is at least $C^1$ and 
$div V = |F|^{n-2} div F +(n-2)|F|^{n-2}\left<\nabla_XF, X\right>$ 
so that the last inequality 
reads
\[
|F|^{n-2}Ricc(X) - div V \leq 0.
\]
 Integrating the last inequality over $M$
we obtain 
\[
 \int |F|^{n-2}Ricc(X) \leq 0.
\]
Since $Ricc (X) \geq 0$ we conclude that at points where $F \neq 0$
the Ricci curvature $Ricc(X) =0$, and hence $F$ is parallel to $E$.
By Proposition~\ref{pnp} the vector field $F$ must be a constant multiple
of $E$.
\end{Proof}$\square$

Proposition~\ref{fsnp} can be applied to the Cartesian product
$\mathbb S^1\times \mathbb S^{n-1}$, and we recover the result 
from \cite{W3}.

\section{Stretched non-positivity}

For a given vector field $E$
we consider the family of Weyl connections defined by vector fields
$\gamma E$, $\gamma > 0$.  For different choices of the background
metric we get different families of Weyl connections. 

\begin{definition}\rm\label{snpwc}
A vector field $E$ is called 
{\it stretched non-positive} (SNP) if there is $\gamma_0 \geq 0$ such 
that the Weyl connections defined by 
the fields $\gamma E$
are non-positive for all $\gamma \geq \gamma_0$.

A Weyl connection on a compact Riemannian manifold $M$ is 
called {\it stretched non-positive} (SNP) if 
its Gauduchon representation $E$ is SNP.
\end{definition}
We have established in Proposition~\ref{snp}  that a parallel vector field
is  SNP.

Our first task is to describe properties of vector fields
$E$, and Weyl connections, with the SNP property.  

\begin{proposition}\label{mrsnp}
If a Weyl connection defined by a divergence free vector field $E$ on
a compact manifold is SNP then 

(W1) $K(\Pi) \leq 0$ for every plane $\Pi \ni E$, where $E \neq 0$;

(W2) $\left<\nabla_YE, Y\right> = 0$ for every $Y\perp E$,
where $E \neq 0$;

(W3) $E$ has constant length on its integral curves.  
\end{proposition}

This Proposition gives us the motivation to restrict our attention 
to divergence free vector fields with constant length, for example 
$E^2 =1$. For such unit vector fields we have

\begin{proposition}\label{mtsnp}
If a unit divergence free vector field $E$
is SNP then  the following conditions are satisfied

(W1) $K(\Pi) \leq 0$ for every plane $\Pi \ni E$,

(W2) $\left<\nabla_YE, Y\right> = 0$ for every $Y\perp E$,

(W4) $\left<\nabla_{E}E, Y_1 \right>^2 \leq 4K(E,Y_2)$ for every
orthonormal frame $(E,Y_1,Y_2)$.

If a unit divergence free  vector field $E$ on a compact manifold
satisfies the conditions (W1), (W2) and

(W5) $\left<\nabla_{E}E, Y_1 \right>^2 <  -4K(E,Y_2)$ for every
orthonormal frame $(E,Y_1,Y_2)$

then the vector field is SNP. 
\end{proposition}
Since the proofs of these two Propositions involve the same 
formulas we will give an amalgamated proof.

\begin{Proof}
First we prove Proposition~\ref{mrsnp}.
Let $X = \frac{E}{|E|}$ where $E\neq 0$.
For unit  mutually orthogonal vectors $X,Y_1,Y_2$ we
consider the plane $\Pi$ 
spanned by $aX+bY_1, a^2+b^2 =1,$ and $Y_2$.
Using Lemma 2.5 we get
for the Weyl connection defined by $\gamma E$
\begin{equation}\rm\label{qf}
\begin{aligned}
&\widehat{K}(\Pi) =
a^2\left(\left<R(Y_2, X)X,Y_2\right> -\gamma(\left<\nabla_{X}E, X \right> +
\left<\nabla_{Y_2}E, Y_2 \right>) \right)+\\
& b^2\left(\left<R(Y_2, Y_1)Y_1,Y_2\right>-
\gamma(\left<\nabla_{Y_1}E, Y_1 \right> + \left<\nabla_{Y_2}E, Y_2 \right>)
-\gamma^2E^2\right)\\+& 
ab\left(2\left<R(Y_2, X)Y_1,Y_2\right> -\gamma
(\left<\nabla_{X}E, Y_1 \right>  + \left<\nabla_{Y_{1}}E, X \right>)
\right).
\end{aligned}
\end{equation}

If the quadratic form (\ref{qf}) is negative semi-definite 
then the coefficient with $a^2$ must be non-positive for all 
$\gamma \geq \gamma_0$. 
It follows that
\[
\left<\nabla_{X}E, X \right>  + \left<\nabla_{Y_2}E, Y_2 \right> \geq 0,
\]
for every $Y_2$ orthogonal to $E$.
Summing the last inequality over an orthonormal  basis in the
subspace orthogonal to $E$, we get   
\[
(n-2) \left<\nabla_{X}E, X \right>  + div E \geq 0.
\]

Since $div E = 0$ we conclude that $ \left<\nabla_{X}E, X \right>  \geq 0$.
This inequality holds where $E \neq 0$, but it implies the inequality
$ \left<\nabla_{E}E, E \right>  \geq 0$ which holds everywhere. 
However the last function is 
the derivative of a smooth function 
$\nabla_E\left(\frac{1}{2}E^2\right) \geq 0$.
Since $div E =0$ the integral of the derivative over the whole manifold
must vanish. It follows that  $\left<\nabla_{X}E, X \right> = 0$  
and $\left<\nabla_{Y}E, Y \right> \geq 0$ for any $Y\perp E$.
However since $ div \ E=0$ we get further that
$\left<\nabla_{Y}E, Y \right> = 0$ for any $Y\perp E$.
This gives us (W3) and (W2).

Now the coefficient with $a^2$ is just the sectional curvature  
of the plane through $E$ and $Y_2$, and hence it must be non-positive,
for every $Y_2$ orthogonal to $E$. This gives us (W1).


To prove Proposition~\ref{mtsnp} we consider 
an orthonormal frame $E,Y_1,Y_2$, and a plane  
plane $\Pi$ spanned by $aE+bY_1, a^2+b^2 =1,$ and $Y_2$.
Since we assume the vector field $E$ to have constant length $1$
the formula \ref{qf} simplifies to  
\begin{equation}\rm\label{qrf}
\begin{aligned}
&\widehat{K}(\Pi) =
a^2\left(\left<R(Y_2, X)X,Y_2\right> -\gamma
\left<\nabla_{Y_2}E, Y_2 \right> \right)+\\
& b^2\left(\left<R(Y_2, Y_1)Y_1,Y_2\right>-
\gamma(\left<\nabla_{Y_1}E, Y_1 \right> + \left<\nabla_{Y_2}E, Y_2 \right>)
-\gamma^2\right)\\+& 
ab\left(2\left<R(Y_2, X)Y_1,Y_2\right> -\gamma
\left<\nabla_{X}E, Y_1 \right>
\right).
\end{aligned}
\end{equation}

If the quadratic form (\ref{qrf}) is negative semi-definite 
for all $\gamma \geq \gamma_0$, then in particular the coefficient 
with $a^2$ must be non-positive.
It follows that  
\[
\left<\nabla_{Y_2}E, Y_2 \right> \geq 0,
\]
for every $Y_2$ orthogonal to $E$. 
Summing the last inequality over an orthonormal  basis in the
subspace orthogonal to $E$, we get   
$div \ E \geq 0$.
Since $E$ is actually a divergence free vector field
we conclude that the terms must be  zero before the summation. 
The condition (W2) follows.

Taking this into account we obtain the following expression 
for the quadratic form (\ref{qrf})
\[
\begin{aligned}
&\widehat{K}(\Pi) =
a^2K(E,Y_2) +  b^2\left(K(Y_1,Y_2) -\gamma^2\right)+\\
&ab\left(2\left<R(Y_2, E)Y_1,Y_2\right> -\gamma
\left<\nabla_{E}E, Y_1 \right> \right).
\end{aligned}
\]

Now the coefficient with $a^2$ is just the sectional curvature   
of the plane through $E$ and $Y_2$, and hence it must be non-positive,
for every $Y_2$ orthogonal to $E$. This gives us (W1).

Taking into account that the discriminant of the form must be 
non-positive for all $\gamma \geq \gamma_0$ we obtain the 
condition (W4).

Conversely if the condition (W5) is satisfied then the discriminant is positive 
for all  $\gamma \geq \gamma_0$. (Note that this calculation
is done pointwise and hence the discriminant is considered on the 
compact grassmanian of tangent planes.) 
Together with conditions (W1) and (W2), this establishes that
the unit vector field $E$ is  SNP.      
\end{Proof}$\square$

It follows from the formula (\ref{qrf}) that 
\begin{proposition}\label{ow} 
If a unit vector field $E$ 
satisfies 
$\left<\nabla_YE, Y\right> > 0$ for every $Y\perp E$
then the field $E$ is SNP. 
\end{proposition}

\section{Homogeneous Weyl connections on groups} 

Let $G$ be a Lie group and $\mathfrak g$ its Lie algebra. 
We endow the group with a left invariant metric. We will denote
the respective scalar product by $< \cdot, \cdot >$.

We have  the following
characterization of left invariant vector 
fields which are parallel or Killing. 
\begin{proposition}\label{gpp}

A left invariant vector field $E$ is a Killing vector field 
if and only if $ad_E$ is skew symmetric.  

A left invariant vector field $E$ is parallel 
if and only if  $ad_E$ is skew-symmetric and 
$E$ is orthogonal to the 
commutator subalgebra $[\mathfrak g, \mathfrak g]$.
\end{proposition}

\begin{Proof}
For left invariant vector fields $E,Y_1,Y_2$
we have 
\[
\left<\nabla_{Y_1}E,Y_2\right> =
\frac{1}{2}\left(\left<\left[Y_1,E\right],Y_2\right> +
\left<\left[Y_2,E\right],Y_1\right> 
+\left<\left[Y_2,Y_1\right],E\right>\right) .
\]
This formula shows that $ad_E$ has the same symmetric part as
the operator $Y\mapsto \nabla_YE$. Hence $E$ is a Killing field
if and only if $ad_E$ is skew-symmetric.

Also the rest of the Proposition follows immediately from this formula.
\end{Proof}$\square$

Note also that a left invariant vector field $E$ is locally a gradient 
of a smooth function (i.e., the form $\varphi = \langle E,\cdot\rangle$
is closed) if and only if $E$ is orthogonal to the commutator 
subalgebra $[\mathfrak g, \mathfrak g]$.

For a left invariant vector field $E$
we consider the Weyl connection defined 
by the field. 

We will study left invariant vector fields $E\in \mathfrak g$
which are SNP, without requiring that $div \ E =0$.
Let us remind the reader that this  is the property of 
the Gaussian thermostat, rather then the Weyl connection,
as discussed in the Introduction.  
However if the field can be factored on a compact 
manifold by a discrete subgroup of isometries 
then by necessity $div \ E =0$, and the respective Weyl
connection is SNP.

Let us recall that a Lie
group is called unimodular if all left invariant vector fields 
are divergence free.

First of all let us consider the 3 dimensional case of unimodular groups.
All the metrics on such groups can be parametrized by three parameters
$\lambda_1,\lambda_2,\lambda_3$, so that for an orthonormal basis 
$(e_1,e_2,e_3)$ in $\mathfrak g$ we have 
\[ 
[e_2,e_3] = \lambda_1 e_1, [e_3,e_1] = \lambda_2 e_2, [e_1,e_2] = \lambda_3 e_3.
\]
The derivation of this can be found in \cite{Mi}.

We searched this class of  examples for non-positive 
and SNP homogeneous Weyl connections. Assuming that not all $\lambda's$ are 
zero we found one SNP example and one isolated non-positive example.
Both occur in the special case of solvable groups $\lambda_3 =0$.

The SNP example occurs when $\lambda_1 = \lambda_2$ and $E = e_3$.
In this case the metric is flat and the vector field is parallel,
see Theorem 4.1 \cite{Mi}. Up to a permutation of indices this is 
the only case of an SNP vector field on a 3 dimensional unimodular 
group. It can be verified easily using our main result 
(see Theorem 4.3 below),
which implies that on a unimodular group the only SNP left invariant vector 
fields are parallel. On the other hand we found direct calculations 
prohibitively cumbersome.

The isolated non-positive example occurs for the solvable group (SOL)
when $\lambda _1 = -\lambda_2 = \lambda$ and $E = \pm \lambda e_3$.
The vector field $E$ cannot be perturbed to any other left invariant
field on SOL defining a non-positive Weyl connection.
We will establish this fact in Section 5 in a more general multidimensional 
setting.

Due to cumbersome formulas we were 
unable to establish that, apart from the flat case,
SOL is indeed the only 3 dimensional 
unimodular Lie group with a non-positive homogeneous Weyl connection.
We conjecture that this is so.

Let us further consider all 4 dimensional extensions
$\mathfrak h$ of the above 3 dimensional unimodular Lie algebra 
$\mathfrak g \subset \mathfrak h$.
We want to check the non-positivity properties of the field 
$E = b \in \mathfrak h$, where $b$ is the unit vector orthogonal 
to $\mathfrak g$.
The extension is completely described by the operator $ad_b = L$. 
However the operator is not arbitrary. 
The Lie algebra $\mathfrak h$ is unimodular if and only if 
$L\mathfrak g \subset \mathfrak g$ (or equivalently 
$[\mathfrak h, \mathfrak h] \subset \mathfrak g$).

By direct calculation
we obtain that the Jacobi identity implies the following

\begin{proposition}\label{laex}
For the diagonal operator $\Lambda$  with eigenvalues 

$1,\lambda_1,\lambda_2,\lambda_3$ the composition 
$L\Lambda$ is anti-symmetric.

In particular if all $\lambda$'s are different from zero
then by necessity $L\mathfrak g \subset \mathfrak g$ and
$[\mathfrak h, \mathfrak h] = \mathfrak g$. Hence 
the Lie algebra $\mathfrak h$ must be unimodular. 
\end{proposition}

By the proof of Proposition  4.1 
the vector field  $E = b$ satisfies (W2) 
if and only if $L$ is anti-symmetric on $\mathfrak g$.
We can see that in the case when $\mathfrak h$ is not unimodular,
$E=b$ is never a Killing vector field.

When $\mathfrak h$ is unimodular we get that $E =b$
satisfies (W2) if and only if it is parallel. If we assume 
$L\neq 0$ then it happens only for special cases:
$\lambda_1 = \lambda_2 = \lambda_3$ and $L$ is an arbitrary
anti-symmetric operator on $\mathfrak g$, or (up to permutation
of indices)
$\lambda_1 =\lambda_2 \neq \lambda_3$ and 
$Le_1 = ae_2, Le_2 = -ae_1, Le_3 = 0$.
We can thus see that parallel vector fields are found  
on nonabelian groups with non-flat left invariant metrics.

All the examples of SNP left invariant vector fields  
described above  in the 3 and 4 dimensional cases
of unimodular groups are parallel.
This is in agreement with our general result

\begin{theorem}\label{gsnp}
For a unimodular Lie group $G$
if a left invariant vector field 
$E \in \mathfrak g$ satisfies properties (W1) and (W2) 
of Proposition 3.1 then it is parallel.
\end{theorem}
It follows that on a unimodular group 
a left invariant field $E\neq 0$ is SNP 
if and only if it is parallel.
We will prove this theorem in the more general setting of 
a homogeneous Riemannian manifold in Section 7.
Before that we will introduce an interesting class
of multidimensional examples. They will show in particular that 
the assumption of unimodularity in Theorem 4.3 is essential.
There are SNP vector fields on non-unimodular groups 
which are not Killing.

\section{Non-positive Weyl connections on extensions of abelian groups}

Let us consider an $n+1$ dimensional extension $\mathfrak s$ 
of an abelian $n$ dimensional 
Lie-algebra $\tilde{\mathfrak s}$ defined by an operator 
$L: \tilde{\mathfrak s}\to  \tilde{\mathfrak s} $. We introduce 
a scalar product into $\mathfrak s$, and let  
$b \in \mathfrak s$ be a unit vector orthogonal to $\tilde{\mathfrak s}$. 
We put $ad_b(u) = [b,u] = Lu, u \in \tilde{\mathfrak s}$,
and we extend naturally the operator $L$ to the whole Lie algebra
$\mathfrak s$. The Lie algebra $\mathfrak s$ is unimodular 
if and only if the trace of $L$ vanishes.

We restrict our attention to the case of the symmetric operator
$L$ with the basis of eigenvectors $e_i, i=0,1, \dots, n$
and respective eigenvalues $\mu_i, i = 0,1,\dots, n$, where 
$e_0 = b$ and $\mu_0 =0$,  

The formulas for the covariant derivative of left 
invariant fields in this case can be found in \cite{Mi}.
We can summarize them as
\begin{equation}\label{cdig}
\begin{aligned}
&\nabla_bb=0, \nabla_bu = 0, \ \ \text{for} \ \ u \in \tilde{\mathfrak s}, \\
&\nabla_u b = -Lu, \nabla_uv = 
\langle Lu,v\rangle b, \ \ \text{for} \ \  u,v \in \tilde{\mathfrak s}
\end{aligned}
\end{equation}

For an orthonormal frame $(X,Y)$ in $\mathfrak s$, 
$X = \sum_{i=0}^{n}x_ie_i, Y =  \sum_{i=0}^{n}y_ie_i$,
we get the Riemannian sectional curvature $K(X,Y)$ 
\begin{equation}
\begin{aligned}
&K(X,Y) = -x_0^2\langle LY,LY\rangle  -y_0^2\langle LX,LX\rangle
+ 2x_0y_0\langle LX,LY\rangle  \\ 
&-\langle LX,X\rangle \langle LY,Y\rangle + \langle LX,Y\rangle^2=\\
&-\sum_{i=1}^{n}\mu_i^2(x_0y_i-y_0x_i)^2 
-\sum_{j < k}^{n}\mu_j\mu_k(x_jy_k-y_jx_k)^2.
\end{aligned}
\end{equation}
         
Further for a vector field $E = \gamma \tilde E= 
\gamma\sum_{i=0}^{n}\eta_ie_i,
\gamma > 0, |\tilde E | =1,$
we have the following expression for the Weylian
sectional curvature $\widehat K(X,Y)$

\begin{equation}
\begin{aligned}
&\widehat K(X,Y)- K(X,Y) =\\ 
&\gamma\left(\eta_0\left(\langle LX,X \rangle 
+\langle LY,Y \rangle\right)-
x_0\langle LX,\tilde E \rangle -y_0\langle LY,\tilde E \rangle
\right)\\
&- \gamma^2 \left(1 - \langle X,\tilde E \rangle^2
-\langle Y,\tilde E \rangle^2\right)
\end{aligned}
\end{equation}

Equipped with these formulas we search for non-positive 
homogeneous Weyl connections on the Lie group with the Lie
algebra $\mathfrak s$.

Let the eigenvalues of $L$ be ordered as  
$\mu_1 \leq \mu_2 \leq \dots \leq \mu_n$ and 
let us assume that all these eigenvalues are different from zero
(apart from $\mu_0 =0$).

\begin{theorem}
The Weyl connection defined by $E = \alpha b$
is non-positive if and only if

(1) for $\mu_1 < 0 < \mu_n$:

$\alpha = \mu_1$ and all the negative
eigenvalues are equal to  $\mu_1$, or 
$\alpha = \mu_n$ and all the positive
eigenvalues are equal to  $\mu_n$.

(2) for $\mu_1 > 0$ ($\mu_n < 0$):

$\alpha \leq \mu_1$ ($\alpha \geq \mu_n$).

\end{theorem}
\begin{Proof}
The sectional curvature can be viewed as a restriction
of a quadratic form on the exterior square $\mathfrak s \wedge \mathfrak s$.

We will take advantage of the fact that in our special case $E = \alpha b$
the quadratic form  defined by $\widehat K$ is ``diagonal''. It follows from 
the following 

\begin{lemma}
For two unit orthogonal vectors $X =\sum_{i=0}^nx_ie_i$
and $Y =\sum_{i=0}^ny_ie_i$
we have for $k=0,1,\dots, n$
\[
x_k^2+y_k^2 = \sum_{j\neq k} m_{kj}^2
\]
where $m_{ij} = x_iy_j-x_jy_i$.
\end{lemma}
\begin{Proof}
We will prove the formula for $k=0$.
Let $\Pi$ be the orthogonal projection onto the 
subspace orthogonal to $e_0$. We compare 
the length of $X\wedge Y$ and $\Pi X \wedge \Pi Y$.
We have $\Pi X \wedge \Pi Y = \sum_{0 < i <j} m_{ij} e_i\wedge e_j$
so $|\Pi X \wedge \Pi Y|^2 =\sum_{0 < i <j} m_{ij}^2$.
On the other hand we have for any two vectors $U,V$
$|U\wedge V|^2 = U^2V^2 -\langle U,V\rangle^2$.
Since the unit vectors $X,Y$ are orthogonal we have 
$\langle \Pi X, \Pi Y \rangle = -x_0y_0$ and 
it follows that 
\[
\begin{aligned}
&|\Pi X \wedge \Pi Y|^2 = |\Pi X|^2|\Pi Y|^2 - \langle \Pi X, \Pi Y \rangle^2
\\&= (1-x_0^2)(1-y_0^2) - x_0^2y_0^2 = 1 - x_0^2  -y_0^2 
\end{aligned}
\]
Since $1=|X\wedge Y|^2 = \sum_{0<i <j} m_{ij}^2 + \sum_{s} m_{0s}^2$
we obtain the desired formula.
\end{Proof}$\square$

Using this lemma we will express $\widehat K(X,Y)$ as a sum of 
squares of $m_{ij}, i,j = 0,1,\dots,n$. The coefficient with 
$m_{0k}^2, k>0,$ is equal to $\mu_k(\alpha - \mu_k)$.
The coefficient  with $m_{ij}, 0 < i < j,$ is equal 
to $-(\alpha -\mu_i)(\alpha - \mu_j)$. The non-positivity of the Weyl
connection is equivalent to the non-positivity of all these coefficients.
We obtain readily the desired conclusions. 
\end{Proof}$\square$

\begin{theorem}
If the  Weyl connection defined by a left-invariant vector  
field $E = \sum_{i=0}^{n}\alpha_ie_i,\in \mathfrak s$
is non-positive  then $\alpha_0\notin (\mu_1,\mu_n)$.

If further $\mu_n > 0$ ($\mu_1 < 0$)  and $\alpha_0 \geq \mu_n$
($\alpha_0 \leq \mu_1$) then 
all of the positive (negative) eigenvalues of $L$ are equal to $\mu_n$
($\mu_1)$, 
and $E = \mu_ne_0$ ($E = \mu_1e_0$).
\end{theorem}
\begin{Proof}
The proof consists of four steps.
We put 
$E = \gamma \tilde E= 
\gamma\sum_{i=0}^{n}\eta_ie_i,
\gamma > 0, |\tilde E | =1$, so that $\alpha_0 = \gamma\eta_0$.

Step 1. $\alpha_0 \notin (\mu_1,\mu_n)$

To prove it we will be choosing an orthonormal frame 
$(X,Y)$ in $\tilde{\mathfrak s}$
such that $E-\alpha_0 e_0\in \tilde{\mathfrak s}$ belongs to the
plane spanned by $(X,Y)$ and $\langle LX,Y\rangle$ = 0. 
The Weylian sectional curvature 
in the direction of this plane is by (9) and (10) equal to 
\[
\begin{aligned}
&\widehat K(X,Y) = -(\gamma \eta_0)^2 +
\gamma\eta_0\left(\langle LX,X\rangle + \langle LY,Y \rangle \right) 
+\langle LX,X\rangle \langle LY,Y \rangle\\
& = - \left(\gamma\eta_0 - \langle LX,X\rangle\right)
\left(\gamma\eta_0 - \langle LY,Y \rangle \right).
\end{aligned}
\]

We conclude that if the Weyl connection is non-positive then
$\gamma\eta_0 \leq \langle LX,X \rangle$ or 
$\gamma\eta_0 \geq \langle LY,Y \rangle$. In other words
$\gamma\eta_0$ does not belong to the open interval 
$(\langle LX,X \rangle, \langle LY,Y \rangle)$.

Assuming $\sum_{j=2}^n\eta_je_j\neq 0$, 
we choose $X = e_1$ and 
$Y= \sigma_1\sum_{j=2}^n\eta_je_j$, where 
$\sigma_1$ is the normalizing factor
$\sigma_1^{-2} = 1 -\eta_0^2-\eta_1^2$.  

We get that $\gamma_0\eta_0$ does not belong to the 
interval $(\mu_1,b)$ where $b$
is a weighted average of the eigenvalues $\mu_{2},\dots,\mu_n$.
More precisely 
$b = \sigma_1^2\left(w + \eta_n^2\mu_n\right)$ 
where $w=\sum_{i=2}^{n-1}\eta_i^2\mu_i$.

Similarly, if $\sum_{j=1}^{n-1}\eta_je_j\neq 0$,
we choose the orthonormal frame $(X,Y)$, where 
$X = e_n$ and $Y =\sigma_n\sum_{j=1}^{n-1}\eta_je_j,
\sigma_n^{-2}= 1-\eta_0^2-\eta_n^2$.
We arrive at the conclusion that  $\gamma\eta_0$ does not belong
to the open interval $(a,\mu_n)$ where
 $a = \sigma_n^2\left(\eta_1^2\mu_1 + w\right)$.
It can be checked that $a< b$. 
It follows that
$\gamma\eta_0$ does not belong
to the open interval $(\mu_1,\mu_n)$.

It remains to analyze the remaining cases of 
$\sum_{j=2}^n\eta_je_j = 0$, 
and/or $\sum_{j=1}^{n-1}\eta_je_j = 0$.
In both cases we choose $X =e_1, Y=e_n$ to arrive 
at the same conclusion.

Step 2. If $\mu_1 < 0$ ($\mu_n > 0$)  and 
$\alpha_0 \leq \mu_1$ ($\alpha_0 \geq \mu_n$)
then $\alpha_0= \mu_1$ ($\alpha_0= \mu_n$),  
and if $\mu_k > \mu_1$($\mu_k < \mu_n$)
then $\alpha_k =0$.

Let us assume without loss of generality that $\alpha_0 = \gamma\eta_0 \leq \mu_1 < 0$.
We choose an orthonormal frame $(X,Y)$ where $X = e_1$, 
$Y= \sigma\left(\eta_0e_0 + \sum_{i=2}^n\eta_ie_i\right)$ and
$\sigma$ is the normalizing factor $\sigma^{-2} = 1- \eta_1^2$.
We have that $E$ belongs to the plane spanned by $X$ and $Y$. 
The Weylian curvature in the direction of this plane is equal
by (9) and (10) to
\[
\begin{aligned}
&\widehat K(X,Y) = 
-\sigma^2\eta_0^2\mu_1^2 - 
\sigma^2\mu_1\sum_{i=2}^n\mu_i\eta_i^2 + \gamma\eta_0\mu_1 =\\
& -\mu_1^2 +\sigma^2\sum_{i=2}^n\left(\mu_1^2-\mu_1\mu_i \right)\eta_i^2 
+\gamma\eta_0\mu_1.
\end{aligned}
\] 
The non-positivity of the curvature implies that 
\[
\mu_1^2 \leq \gamma\eta_0\mu_1 \leq \mu_1^2 
-\sigma^2\sum_{i=2}^n\left(\mu_1^2-\mu_1\mu_i \right)\eta_i^2\leq \mu_1^2.
\]
It follows that $\sum_{i=2}^n\left(\mu_1^2-\mu_1\mu_i \right)\eta_i^2 =0$
and $\gamma\eta_0 = \mu_1$. 
We have thus achieved Step 2.

Step 3.  If $\alpha_0 = \mu_1 <0$ then $E = \mu_1e_0$.

We have by step 2 that $E = \gamma\sum_{i=0}^k\eta_ie_i$
and $\mu_1 = \mu_2 = \dots = \mu_k < 0$.
We choose  an orthonormal frame $(X,Y)$ where 
$X = x_0e_0 + x_1e_1$ with arbitrary coefficients $x_0,x_1, x_0^2+x_1^2 =1$,
and $Y= \sigma \sum_{j=2}^k\eta_je_j$, where 
$\sigma^{-2} = 1-\eta_0^2 - \eta_1^2$. If $E = \alpha_0e_0 + \alpha_1e_1$
then we take $Y$ to be an arbitrary unit vector orthogonal to all the vectors 
$X$. The Weylian sectional curvature $\widehat K(X,Y)$ can be represented as
a quadratic form in the variables $(x_0,x_1)$, 
$\widehat K(X,Y) = Ax_0^2 +B x_0x_1 + Cx_1^2$.

By direct calculation we  get  $C = 0$ and $B = \mu_1\alpha_1$.
It follows that $\alpha_1 =0$. Since the status of all the components
$\alpha_1,\alpha_2, \dots, \alpha_k$ is the same under present assumptions, we conclude
that all of them must be zero.

Step 4. If $\mu_1 <0$ and $E = \mu_1e_0$  
then all the negative eigenvalues of $L$ are equal to $\mu_1$.

For $\mu_1 < \mu_i <0 $ we choose the orthonormal frame $X = e_0$ and $Y= e_i$. 
We get that $\widehat K(X,Y) = -\mu_i^2 +\mu_1\mu_i > 0 $. This contradiction
ends the proof of Step 4, and the proof of the Theorem.
\end{Proof}$\square$

The last theorem gives us immediately the following 

\begin{corollary}\label{nnm}
If $L$ has both positive and negative eigenvalues and 
the Weyl connection defined by $E$ is non-positive
then $E= \mu_ne_0$ and all the positive eigenvalues are
equal to $\mu_n$, or 
$E= \mu_1e_0$ and all the negative eigenvalues are
equal to $\mu_1$.
\end{corollary}
The last Corollary covers the case of a unimodular 
non-abelian Lie algebra $\mathfrak s$,
since then $\sum_{i=1}^n\mu_i =0$.

Let us assume that  $e^L$ is an automorphism of a discrete subgroup 
$\Gamma_0 < \mathbb R^n$ with $n$ generators. The group $\Gamma > \Gamma_0$
generated by $\Gamma_0$ and $e^L$ is a discrete subgroup of 
the simply connected Lie group $S$ 
with the Lie algebra $\mathfrak s$. The homogeneous space
$M = _{\Gamma}\S$ is then compact and we get a family of examples 
of compact manifolds with only one isolated non-positive Weyl connection. 
The geodesic flows of such connections were recently investigated 
in \cite{W4}. One of the conclusions is the description 
of the Gaussian thermostat defined on $M$ by the left invariant 
vector field $E =b$. As explained in the Introduction we may get different
dynamics for different magnitudes of the velocity $k = v^2$.
In the case of $E=b$ the Gaussian thermostat for $0 < k < \mu_n^{-1}$
is asymptotic to an Anosov
flow (the suspension of the hyperbolic toral automorphism $e^L$). 
For $k > \mu_n^{-1}$ an open dense subset of 
the phase space is foliated by invariant tori carrying 
quasi-periodic motions, while the Anosov subsystem 
on a submanifold of half the dimension of the phase space is still
present. Hence the isolated non-positive Weyl connection defined
by $\mu_nE$ corresponds to the borderline dynamics.

Finally let us analyze the case of all positive,
(or all negative) eigenvalues, 
$0 < \mu_1 \leq \mu_2 \leq \dots \leq \mu_n$. 
In this case the 
Riemannian sectional curvatures are all negative,
and hence the homogeneous 
Weyl connection defined by $E\in \mathfrak s$
is non-positive if only $\gamma = |E|$ is sufficiently
small. We do not attempt to give explicit conditions
for non-positivity in this case. Instead we describe 
all SNP vector fields $E \in \mathfrak s$ under the assumption
that $ div \ E =0$.

For a unit vector field $E = \sum_{i=0}^n\eta_ie_i$
we have $div E = \eta_0\sum_{i=0}^n\mu_i$. Hence 
either  $\sum_{i=0}^n\mu_i = 0$ and the group is unimodular,
or $\sum_{i=0}^n\mu_i \neq 0$ and then  $div \ E = 0$ if and only 
if $E \in \widetilde {\mathfrak s}$ (i.e.,  
if $E$ is orthogonal to $b = e_0$).

For a vector field $E = -e_0$ we get $\langle \nabla_YE,Y\rangle =1$
for every unit vector $Y$ orthogonal to $E$. It follows from 
Proposition~\ref{ow} that all vectors close to $-e_0$ are SNP.

The divergence free vector fields  $E\in \widetilde{\mathfrak s}$
which are SNP are described in 
the following

\begin{proposition}\label{nnm1}
If a unit divergence free vector field 
$E\in \widetilde{\mathfrak s}$ is  
SNP then $E$ is an eigenvector of
$L$ with the eigenvalue $\mu_k$ such that 
$\mu_k \leq 4 \mu_1$.

If  $E $ is an eigenvector of
$L$ with the eigenvalue $\mu_k$ such that 
$0< \mu_k < 4 \mu_1$ then it is  SNP. 
\end{proposition}
\begin{Proof}
We apply directly  Proposition~\ref{mtsnp}.
For a vector $Y=\sum_{i=0}^ny_ie_i$ orthogonal 
to $E$ we must have $\langle \nabla_YE,Y\rangle =0$.
By (\ref{cdig}) we get $\nabla_YE = \langle LY,E \rangle b$,
and finally if $y_0 \neq 0$ we conclude that
$\langle LY,E\rangle =\langle Y,LE\rangle = 0$ for any vector $Y$ orthogonal 
to $E$. It implies that $E$ must be an eigenvector  
of $L$. 

If $E = e_k$ is an eigenvector of $L$ with the eigenvalue 
$\mu_k$ then for $Y_1 = b$ the condition (W3) says that
$\mu_k^2 \leq -4 K(E,Y_2)$ for any vector $Y_2 \in
\tilde{\mathfrak s}$ orthogonal to $E$.
The Riemannian curvature in the plane spanned by $E$ and
$Y_2= \sum_{i\neq k}^ny_ie_i$ is by (9) equal to 
\[
-K(E,Y_2) =y_0^2\mu_k^2 + \sum_{i\neq 0,k}y_i^2\mu_i\mu_k
\]
It is the smallest in absolute value for $Y_2 = e_1$
and $-K(e_k,e_1) = \mu_1\mu_k$. This gives us 
the condition $\mu_k \leq 4 \mu_1$.

Conversely, if $\mu_k <  4 \mu_1$ then the condition
(W5) is satisfied and the proof is complete. 
\end{Proof}$\square$

The above examples show that
a unit divergence free vector field $E$ can have the SNP 
property without being a Killing field.

In particular if $L$ is equal to identity then we get 
the $(n+1)$-dimensional hyperbolic space $\mathcal H^{n+1}$.
If we model it  by the upper half space 
then the Lie group of isometries consists of translations
and dilations, and it can be identified with
$\mathcal H^{n+1}$.  

The unit divergence 
free left invariant vector fields are horizontal fields
which have constant direction in the Euclidean sense, e.g. it has 
the direction of the first axis. Clearly they do not have 
constant length in the Euclidean sense and they are not Killing
fields. Any such vector field defines a homogeneous Weyl connection which
is SNP. 
Let us note that this example cannot be factored 
onto a compact manifold. 
We do not know an example of a 
divergence free vector field  on a compact manifold 
with the SNP property which is not parallel. 

Additionally in the hyperbolic space we have that
for a vector field $E = \sum_{i=0}^n\alpha_i e_i$
we have $\langle \nabla_YE,Y\rangle = -\alpha_0$
for every unit vector $Y$ orthogonal to $E$.

Hence by Proposition~\ref{ow} we have in the hyperbolic space 
that all vector fields $E$ with $\langle E,b \rangle < 0$
are SNP.
These fields are clearly  not Killing vector fields.

\section{Homogeneous Riemannian spaces}

Let $\pi: G \to G/H = M$ be a Riemannian homogeneous space.
We assume that the Lie algebra $\mathfrak g = \mathfrak p \oplus \mathfrak h$,
where $\mathfrak h$ is the Lie algebra of the compact subgroup $H < G$, 
and the splitting is $ Ad(H)$ invariant. 

Let $e$ denote the unit element in $G$ and let $\pi(e) = o$. 
We can identify the tangent space $T_oM$ with $\mathfrak p$.
However let us note that at any other point in $M$ we do not
have in general a canonical identification of the tangent space with
$\mathfrak p$.
We choose a scalar product in $\mathfrak p$ which is $Ad(H)$ invariant.
We spread the scalar product to other tangent spaces of $M$ by the 
action of the Lie group $G$ on the manifold $M= G/H$.
This  gives us {\it  a homogeneous Riemannian metric on $M$}.

The left-invariant  vector fields on $G$ do not project onto $M$.
However the right-invariant vector fields do. Indeed, the action
by left translations
of a one parameter subgroup $g^u, u \in \mathbb R,$ 
projects onto $M$ as a group 
of isometries. Its velocity field is hence a Killing vector field on $M$.

These Killing fields are projections under $\pi:G \to M =G/H $ of the right-invariant 
vector fields on the group $G$. 
Choosing an element from the Lie algebra $\mathfrak g$ determines 
uniquely the whole Killing field on $M$.  
The Killing field  vanishes usually at some points
away from the point $o$.

For Killing vector fields $X,Y,Z$ on $M$ we have the following

\[
\left\langle \nabla_X Y , Z \right\rangle = 
\frac{1}{2} \left( 
\left\langle \left[ X, Y \right] , Z \right\rangle - 
\left\langle \left[ Z, X \right] , Y \right\rangle - 
\left\langle \left[ Z, Y \right] , X \right\rangle \right)
\]

It is a completely general formula of Riemannian geometry.
To apply this formula in our special setting we need to exercise 
some caution since the covariant derivative of a Killing vector field
is not necessarily Killing. We will consider the last formula 
where  $X,Y$ and $Z$ are projections on $M$ of 
right-invariant vector fields on $G$.  
Such vector fields project as Killing vector fields on 
$M$. The commutator of 
right-invariant vector fields projects to $M$ as the commutator 
of respective projections. It allows us to obtain an explicit 
formula for the covariant derivative, at least at the point 
$o$.

We assume $Z\in \mathfrak p$, however we allow arbitrary 
right-invariant vector fields $X$ and $Y$.

We introduce the tensor $U(X,Y)$ in  
$\mathfrak g$ with values in $\mathfrak p$ by

\begin{equation}\label{wnu}
\left\langle U(X,Y) , Z \right\rangle = 
\frac{1}{2} \left( 
\left\langle \left[ Z,X  \right]_{\mathfrak p} , Y \right\rangle + 
\left\langle \left[ Z, Y \right]_{\mathfrak p} , X \right\rangle \right),
\end{equation}
where the subscript $\mathfrak p$ indicates the projection  to $\mathfrak p$.
At the point $o$ we get for two Killing vector fields on $M$, 
obtained as projections of 
right-invariant vector fields $X,Y$ on $G$ 
\begin{equation}\label{cdih}
\left(\nabla_X Y\right)(o)  = 
-\frac{1}{2}\left[ X, Y \right]_{\mathfrak p}  + U(X,Y). 
\end{equation}

We will need a formula for the sectional curvature on $M$. It can be written 
in terms of a tensor on $\mathfrak p$. It can be found in \cite{Be}, with
a complete derivation. However we will need it in a somewhat different
form, and we give an independent proof. Note the apparent asymmetry
of the following formula.  

\begin{theorem}
The sectional curvature $K = K(v,J)$ of $M$ 
in the direction of the plane spanned by  
an orthonormal frame $v,J \in \mathfrak p$ is equal to
$$
\aligned
K =
&\left(\frac{1}{2}[J,v]_{\mathfrak p}+U(v,J) \right)^2 \\
&-\left\langle \left[ \left[  J, v \right], 
v \right]_{\mathfrak p} , J \right\rangle 
+\left\langle \left[   J, U(v,v) \right]_{\mathfrak p} , J \right\rangle 
-\left\langle \left[   J, v \right]_{\mathfrak p} , 
\left[   J, v \right]_{\mathfrak p} \right\rangle
\endaligned
$$
\end{theorem} 
\begin{Proof}
We will use the Jacobi fields of the Riemannian geometry.

Let $\gamma_t  \in G, t\in (-\epsilon, \epsilon), $ be
 a lift of a geodesic from $M$ to $G$, 
where $\gamma_0$ is the unit element in $G$, and $v(t) = \frac{d}{dt} \gamma_t$ is the 
velocity field along $\gamma_t$.

The basic observation is that the restriction of any Killing vector field  
to the geodesic $\gamma_t$ is a Jacobi field.

Let $J\in \mathfrak p,$ be a right-invariant vector field in $G$. 
Its projection to $M$ is a Killing field, and its restriction to the 
geodesic is a Jacobi field. We will denote the resulting Jacobi field
also as $J$. If $J\in \mathfrak p$ is a unit vector field orthogonal
to $v = v(0) \in \mathfrak p$ then we get from the Jacobi equations
the following general formula for the sectional curvature $K=K(v,J)$ 
$$
K= \left(\nabla_vJ\right)^2 - \frac{1}{2}\frac{d^2}{dt^2} J^2.
$$

To get the first term we can use the formula (\ref{cdih})
\[
\nabla_vJ = -\frac{1}{2}[v,J]_{\mathfrak p}+U(v,J).
\]
However since it only holds at the point $o$ we cannot use it directly
to calculate $\frac{d^2}{dt^2} J^2$. 
This is the main stumbling block in our derivation of the formula.

We move the Jacobi vector field $J(t) = J(\gamma_t)$ to $\mathfrak g$
by left translations and denote the resulting function by 
$\widetilde J: (-\epsilon, \epsilon) \to \mathfrak g$. We do not claim that
$\widetilde J$ depends only on the geodesic in $M$, 
it actually depends on the particular lift $\gamma_t$ of the geodesic. 
However we can safely claim that $\left(J(t)\right)^2 = 
\left(\widetilde J(t)\right)^2$ because the left translations on $G$
project as isometries of $M$. Our task now is to differentiate the function
$\left(\widetilde J\right)^2$.

We will thus differentiate $\widetilde J(t)$ twice with respect to the real 
parameter $t$. In the calculation of $\left(\widetilde J(t)\right)^2$
only the projection onto $\mathfrak p$ matters, but we can differentiate
first and project later, since the differentiation is done in the linear
space $\mathfrak g$.

It is well known how to take the first derivative 
$$
\frac{d}{dt} \widetilde J  = \left[\widetilde J, \widetilde v\right],
$$
where $\widetilde v (t) \in \mathfrak g$ is the left translation 
of the velocity vector field $v(t)$ along the geodesic $\gamma_t$
to the Lie algebra $\mathfrak g$ of $G$. It is a general formula valid in a Lie group,
and it holds in the whole interval $(-\epsilon, \epsilon)$.
Now we differentiate the first derivative only at $t=0$. To achieve that
we need yet to find $\frac{d}{dt} \widetilde v$ at $t=0$.

Let us consider the velocity vector field $v(t)$ as the restriction to 
$\gamma_t$ of a time dependent right-invariant vector field
$\widehat v(t)\in \mathfrak g$. Namely we put 
$\widehat v = Ad_{\gamma_t}\widetilde v$. We have clearly
$$
\frac{d}{dt}_{|t=0} \widehat v = \frac{d}{dt}_{|t=0} \widetilde v.
$$

At the same time using the equations for geodesics and the formula (\ref{cdih}) 
we obtain at the point $o$
$$
0 = \nabla_vv = \frac{d}{dt}_{|t=0} \widehat v + U(v,v).  
$$
Note that we do not claim here that the right hand side 
($\frac{d}{dt}_{|t=0} \widehat v$) belongs to $\mathfrak p$.
The equality is to be understood as the equality of projections
to the subspace $\mathfrak p$.

We obtain thus
$$
\frac{d}{dt}_{|t=0} \widetilde v = -U(v,v) + r,
$$
for some $r \in \mathfrak h$ 
Since the adjoint action of $H$ on $\mathfrak p$ is
by isometries, we have that for any $r \in \mathfrak h$
the operator $ad_r: \mathfrak p \to \mathfrak p$
is skew-symmetric. Since at $t=0$ we have 
$\widetilde J = J\in \mathfrak p$ then it follows that
$$
\langle[\widetilde J, r],\widetilde J\rangle = 
\langle [J,r], J\rangle =0.
$$

Finally we have at $t=0$
$$
\aligned
&K(v,J)- \left(\frac{1}{2}[J,v]_{\mathfrak p}+U(v,J) \right)^2 = 
-\frac{d}{dt} \left\langle \left[ \widetilde J, \widetilde v \right]_{\mathfrak p} , 
\widetilde J \right\rangle = \\
&- 
\langle \left[ \left[  J, v \right], v \right]_{\mathfrak p} ,
J \rangle +
\langle \left[J, U(v,v) \right]_{\mathfrak p} , 
 J \rangle -\langle \left[ J, v \right]_{\mathfrak p} , 
\left[   J, v \right]_{\mathfrak p} \rangle,
\endaligned
$$
and the desired formula follows.
\end{Proof}$\square$

\section{Homogeneous Weyl connections on  homogeneous spaces}

To introduce homogeneous Weyl connections on the Riemannian homogeneous 
space $G/H$ of Section 6 we consider the subspace
\[
\mathfrak p_0 = \{X \in \mathfrak p| [Y,X] = 0 \ \ \forall Y\in \mathfrak h\},
\]
on which the action of $Ad(H)$ is trivial.

We choose a left invariant vector field $E$ defined by an element 
from $\mathfrak p_0$.
It projects onto $M = G/H$ as a $G$-invariant vector field. The subspace
$\mathfrak p_0$ may be equal to zero, however 
for $\mathfrak h = 0$ all of $\mathfrak g$ is in $\mathfrak p_0$. 
The last case is that of a left invariant metric on $G$, with 
a left invariant vector field $E$. 

The Weyl connection defined on $M = G/H$  by the projection of 
a left-invariant vector field $E \in \mathfrak p_0$ is clearly 
preserved by the isometric action of the group $G$ on $M$.

Our main result says

\begin{theorem}
If the Lie group $G$ is unimodular and a $G$-invariant vector 
field on $M= G/H$, defined by $E \in \mathfrak p_0$, is SNP
then the vector field is parallel.
\end{theorem}

We will prove a slightly more general fact.

\begin{theorem}
If the Lie group $G$ is unimodular and a $G$-invariant vector 
field on $M= G/H$, defined by $E \in \mathfrak p_0$, 
satisfies properties (W1) and (W2) then the vector field is parallel.
\end{theorem}
\begin{Proof}
Without loss of generality we can assume $E^2 = 1$

We need a formula for $\nabla_YE$, where $Y\in \mathfrak p$ is orthogonal 
to $E$. The field $E$ can be naturally considered as a left-invariant
vector field 
in $G$ which projects onto $M = G/H$. We can extend $Y$ to a right-invariant
vector field in $G$ which also projects onto $M$ (as a Killing vector 
field). Since any left-invariant vector field commutes with any right-invariant
vector field the  projections  of $E$ and $Y$ on $M$ also commute. 
We conclude that $\nabla_YE = \nabla_EY$ and $Y$ is a Killing vector
field. At the point $o$ we can use the formula (\ref{cdih})
to get
$$
\nabla_EY = \frac{1}{2}[Y,E]_{\mathfrak p} +U(E,Y).
$$
Moreover by the definition of $U$ we get
$$
\langle U(E,Y),Y\rangle = \frac{1}{2}\langle [Y,E]_{\mathfrak p},Y\rangle.
$$

We conclude that 
$$
\langle \nabla_YE,Y\rangle = \langle [Y,E]_{\mathfrak p},Y\rangle.
$$

By the last formula   the condition (W2) translates into
\[ 
\langle [E,Y]_{\mathfrak p}, Y \rangle =0
\] 
for any field $Y\in \mathfrak p$ orthogonal to $E$.

Hence the condition (W2) can be rephrased in the following way. 
Let $A:\mathfrak p \to \mathfrak p$
be the operator defined as the composition of the operator 
$ad_E = [E,\cdot]$ and the orthogonal projection 
onto the orthogonal complement
of $E$. The property (W2) is now equivalent to $A$ being skew-symmetric.
We have 
\[[E,Y]_{\mathfrak p} = A(Y) + \sigma(Y)E
\]
for some linear functional 
$\sigma : \mathfrak p \to \mathbb R$.

We invoke now the formula for the sectional curvature 
of $M$ at the point $o$ in the plane  spanned by 
the unit $G$-invariant vector field $v =E$ and 
the Killing vector field $J= Y \in \mathfrak p$, $Y$ orthogonal to $E$.
\[
K(E,Y) =
(\nabla_YE)^2 -
\left<\left[E,\left[E,Y\right]\right]_{\mathfrak p},Y\right> -
\left<\left[U(E,E),Y \right]_{\mathfrak p},Y\right> -[E,Y]_{\mathfrak p}^2
\]

Under the condition (W2) we have further, using the fact that $ad_{\mathfrak h}$
vanishes on $E$,
\[
\begin{aligned}
&\left<\left[E,\left[E,Y\right]\right]_{\mathfrak p},Y\right> 
=\left<\left[E,\left[E,Y\right]_{\mathfrak p}\right]_{\mathfrak p},Y\right> 
=\left<\left[E,A(Y)\right]_{\mathfrak p},Y\right>= \\
&=-\left<\left[E,Y\right]_{\mathfrak p},A(Y)\right> 
=-\left<\left[E,Y\right]_{\mathfrak p},\left[E,Y\right]_{\mathfrak p}\right> +
\sigma(Y)^2. 
\end{aligned}
\]

Now the formula for the curvature reads
\begin{equation}
\label{fke}
K(E,Y) =
(\nabla_YE)^2 -\sigma(Y)^2 -
\langle\left[U(E,E),Y\right]_{\mathfrak p},Y\rangle.
\end{equation}

Let us choose in $\mathfrak p$ an orthonormal basis
$E, Y_1, Y_2, \dots Y_{n-1}$.
 
We calculate the Ricci curvature in the direction of $E$
by summing the sectional curvatures over the chosen 
orthogonal basis. It follows from the condition (W1) that
$Ricc(E)\leq 0$.

\begin{equation}
\label{Ric}
\begin{aligned}
&0 \geq Ricc(E) = \sum_{i=1}^{n-1}K(E,Y_i)=\\
&=\sum_{i=1}^{n-1}(\nabla_{Y_i}E)^2  
-\sum_{i=1}^{n-1}\sigma(Y_i)^2 
-\sum_{i=1}^{n-1}\langle\left[U(E,E), Y_i\right]_{\mathfrak p}, Y_i\rangle.
\end{aligned}
\end{equation}

Invoking the unimodularity of the Lie group we have 
that $ad_{U(E,E)}$ is traceless. Since $ad_{\mathfrak h}$ preserves $\mathfrak p$ 
there is no contribution to the trace from $\mathfrak h$ and we conclude that
\begin{equation}
\label{tB}
tr \ ad_{U(E,E)}= 
\sum_{i=1}^{n-1}\langle\left[U(E,E),Y_i\right]_{\mathfrak p}, Y_i\rangle 
+\langle\left[U(E,E),E\right]_{\mathfrak p}, E\rangle = 0
\end{equation}

Further we observe that in view of (\ref{wnu})   we get
\[
U(E,E) =  - \sum_{i=1}^{n-1}\sigma(Y_i)Y_i \ \ \ \
\text{and}
\ \ \ \ \
\langle\left[U(E,E),E\right]_{\mathfrak p}, E\rangle = 
\sum_{i=1}^{n-1}\sigma(Y_i)^2 
\]

Hence using this and (\ref{tB})  the inequality (\ref{Ric})  becomes
\[
0 \geq Ricc(E) =
\sum_{i=1}^{n-1}(\nabla_{Y_i}E)^2.
\]

We conclude that for every $Y$ orthogonal to $E$
\[
\nabla_YE =0 \ \ \ \text{and} \ \ \ K(E,Y) = 0.
\]

Going back to the equation (\ref{fke}) we obtain for every
$Y \perp E$
\[
0 = K(E,Y) =
-\sigma(Y)^2 -
\langle\left[U(E,E), Y\right]_{\mathfrak p},Y\rangle.
\]
Substituting $Y =U(E,E)$ we obtain 
\[
0 =\sigma(U(E,E)) = -\sum_{i=1}^{n-1}\sigma(Y_i)^2
\]
and hence the functional $\sigma$ vanishes. We conclude
that also $\nabla_EE = U(E,E) = 0$. We can see that 
$\nabla_YE = 0$ for every vector $Y\in \mathfrak p$, i.e., 
the vector field $E$ is parallel on $M$.
\end{Proof}$\square$

\end{document}